\documentclass[a4paper,12pt]{article}

\usepackage[latin1]{inputenc}
\usepackage[francais,english]{babel}
\usepackage{graphicx,url}
\usepackage{latexsym,amssymb,amsmath,amsthm}
\usepackage[T1]{fontenc}

\theoremstyle{definition}

\theoremstyle{remark}

\numberwithin{equation}{section}

\begin{document}

\title{Jakob Bielfeld (1717--1770) and the diffusion of statistical
  concepts in eighteenth century Europe}
\author{Bernard Ycart\footnote{
Laboratoire Jean Kuntzmann,
Universit\'e Grenoble Alpes and CNRS UMR 5224, 
51 rue des Math\'ematiques 38041 Grenoble cedex 9, France
\texttt{Bernard.Ycart@imag.fr}}}
\maketitle
\begin{abstract}
Published between 1760 and 1770, Bielfeld's writings prove that
scholars of the time were acquainted with the concepts of
both political arithmetic and German statistik, long before they
merged into a new discipline at the beginning the following
century. It is argued here that these works may 
have been an important source of diffusion of
statistical concepts at the end of 
the eighteenth century. Bielfeld is now almost completely forgotten,
and the reasons for his lack of fame in posterity are examined.
\end{abstract}
\selectlanguage{francais}
\begin{abstract}
Publi\'es entre 1760 et 1770, les \'ecrits de Bielfeld prouvent que
les concepts de l'arithm\'etique politique et de la statistik 
allemande \'etaient largement r\'epandus parmi les savants, longtemps
avant que les deux ne fusionnent en une nouvelle discipline au d\'ebut
du si\`ecle suivant. On soutient ici que ces publications peuvent
avoir \'et\'e une source importante de diffusion des concepts
statistiques \`a la fin du dix-huiti\`eme si\`ecle. Bielfeld est de
nos jours presque compl\`etement oubli\'e, et les raisons de cet oubli
sont examin\'ees. 
\end{abstract}
\selectlanguage{english}

{\small

}
%
\section{Introduction}
%
In the autumn term of 1921, Karl Pearson (1857-1936) 
started a series of lectures
on the history of statistics. In the introduction, 
he vividly described the origins of the discipline
and its naming.
\begin{quotation}
To this hybrid discipline of statecraft, constitutional history and
description of state constitutions, Gottfried A. Achenwall in 1752, for
the first time as I am aware, introduced the word `Statistik' as the
name of a distinct branch of knowledge. That branch was not concerned
with numbers nor with mathematical theory.

Meanwhile in England there was an entirely different movement.
Captain John Graunt [\ldots], who lived from 1620 to 1674, a clear
century before Achenwall (1719-1772) [\ldots and his] 
friend Sir William Petty [\ldots] founded
the English school of what was called `Political Arithmetic'.

A hundred years and more later comes an extraordinary event. A
Scotsman steals the words `Statistics' and `Statistik' and applies
them to the data and methods of `political arithmetic'. It was
certainly a bold, bare-faced act of robbery which Sir John Sinclair
commited in 1798  [Pearson 1921:2].
\end{quotation}
We now know that viewing English political arithmetic 
and German statistik as ignorant from each other, 
before Sinclair ``boldly'' synthesized them
into a new discipline at the end of the eightenth century, is
too schematic. In the first half of the eighteenth
century, political arithmetic was a declining 
discipline in England
[Studenski 1961:40; Deane 1987];
its revival in the second half of the century was largely 
continental [Todhunter 1865; Hacking 1984]. 
In Germany, even though the G\"ottingen
school of statistik sometimes opposed mathematics 
[John 1883:670; Hacking 1990:24], 
political arithmetic was regarded as one methodological
component in the very 
broad definition of statistik. It was actively developed in particular
by Johan Peter S\"ussmilch (1707-1767) [Heuschling 1845:8; Hacking 1984:113].
John Sinclair's monumental \emph{Statistical Account of Scotland}
`was not a work of statistics in the modern sense but it was one
in a looser version of the old German sense' [Cullen 1975:10]. Neither was it the 
first one of the sort in Great Britain. Translations of German statistik books
were available in English before Sinclair, such as 
[B\"usching 1762; Zimmerman 1787]. Arthur Young (1741-1820), who was
perfectly acquainted with the techniques of political arithmetic, published
his \emph{Six month tour through the North of England} in 1770 
[Young 1770, De Bruyn 2004], and it was
recognized in Germany as a true work of statistik [Meusel 1790:24]. 

By the time Sinclair perpetrated his ``bare-faced act of robbery'',
the concepts, vocabulary, and methods of both political arithmetic
and statistik had long pervaded the whole Europe of Enlightment.
Ample evidence can be found in the books of a popular-science
writer of the time, Jakob Friedrich Bielfeld (Freiherr von)
(1717-1770), who wrote in French and signed ``Monsieur le Baron 
de Bielfeld''. In his \emph{Institutions politiques} 
[Bielfeld 1760a, 1760b, 1772], and \emph{Les premiers 
traits de l'\'erudition universelle} [Bielfeld 1767, 1770a,b], the
following can be found.
\begin{itemize}
\item a history of political arithmetic up to his time,
\item the main concepts of political arithmetic (mortality tables,
insurance, life annuities, tontines),
\item a methodology for data collection, similar to
the one Sinclair would implement 30 years later,
\item a statistik of Europe,
\item the first use of the words \emph{statistique} 
in French, and \emph{statistic} in English. 
\end{itemize}
The first objective of this article is to review Bielfeld's works 
on political arithmetic and statistik. Bielfeld never viewed himself as a 
scholar, nore did he pretend to expose the result of his own research: 
his writings are just a clear and synthetic exposition of 
what was generally understood in his time. 

Nowadays, Bielfeld's name does not appear in authoritative histories
of statistics,  
he is not among the ``leading personalities 
in statistical sciences'' [Johnson \& Kotz 1997], and does not have an entry
 in Burns' encyclopedia of science in the enlightenment [Burns 2003]. 
 Few books cite him
as being the first to use of the term `statistic' 
[Cullen 1975:10; Nalimov 1981:208; Federer 1991:1;  Headrick 2000:68; 
Agarwal 2009:1], whereas many more sources (e.g. [Hald 2003]) cite Pearson 
and/or reproduce Sinclair's own version
[Sinclair 1798:xiii].
As remarked by Reinert [2013], the man himself is by 
now almost completely forgotten. He has an entry in the 
\emph{Allgemeine Deutsche Biographie} [Steffenhagen 1875] 
and in the German version of Wikipedia, 
but is  scarcely found in any encyclopedia outside Germany, from 1850 on. 
However, historians of economy [Stangeland 1904:303, Tribe 1988:82,
Reinert 2009]
and politics [Bazzoli 1990, S\'anchez-Blanco 2003, Cunha 2011] have
long recognized the influence 
of Bielfeld's writings on western European thinking. 
Here is for instance how Cunha describes it.
\begin{quotation}
It is possible that the author who contributed the most, not to the
formulation but to the diffusion of Cameralism, was Jakob Friedrich
von Bielfeld (1717-70). His \emph{Institutions politiques}, originally
published in French in 1760, was reprinted 12 times, in addition to
being translated into German, Russian, Spanish and Italian. This book
became a kind of economics bestseller in its day, containing for
example an extensive typology of the origins of the decadence of states.
[Cunha 2011:68]
\end{quotation}
There is no reason to believe that Bielfeld's influence was smaller for
statistics than it was for politics and economy, and the second objective 
of this article is to assess its importance. 
Evidence from various sources that Bielfeld was
read, copied,
commented, and sometimes hotly debated, will be presented. 
Particular emphasis will be put on Bielfeld's controversy with
S\"ussmilch, which probably fueled Euler's interest in political
arithmetic, entailing important consequences 
for the development of the discipline 
[Todhunter 1865:239-247; Klyve 2013].

Our third objective is an attempt at answering the question of 
fame in posterity: why did Bielfeld's name nearly disappear from the 
history of statistics? Several reasons will be
identified: Bielfeld's mitigated eulogy by Formey at the Berlin academy
[Formey 1770], was an easy argument for his detractors; 
his theatrical works were harshly
criticised in his time, and his unreliable recollections
from Frederick's court in Rheinberg were severely judged by 
nineteenth century historians; finally, his short life
after the publication of his main works, and perhaps his own personality,
may have prevented him from maintaining a dense network 
of correspondents, and ensuring his own publicity.

The article is divided into three sections, dealing with
biography, political arithmetic, and statistik. Each of the 
three sections contains two subsections, the first one
on primary sources, the second one on posterity.  

\section{Bielfeld and his fading fame}
%
\subsection{Who was Bielfeld?}
\begin{figure}[!ht]
\centerline{
\includegraphics[width=10cm]{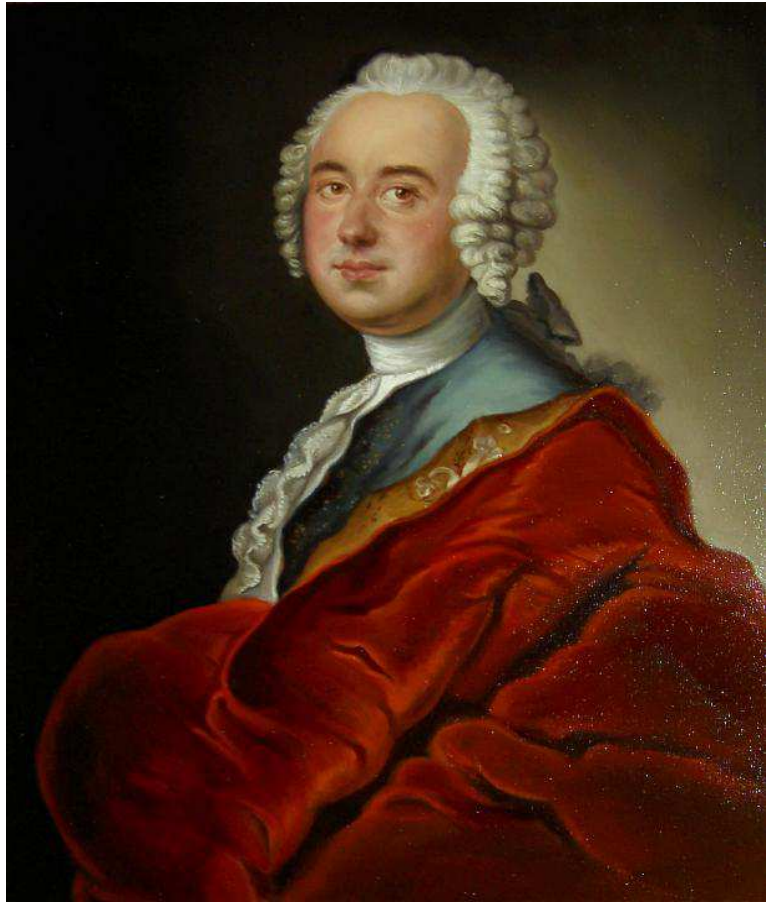}
} 
\caption{Portrait of Jakob Friedrich Bielfeld (Freiherr von)
  (1717-1770).
Source: Guido Sch\"aferhoff (Wikipedia).}
\label{fig:Color}
\end{figure}
Erik Reinert gives the following summary of Bielfeld's life and 
career.
\begin{quotation}
Born into a family of merchants in Hamburg in 1717, Jakob Friedrich von
Bielfeld's life was international already from an early age. In 1732 he
started his university studies in Leyden, Holland. In 1735 he travelled in
Holland, France, and England. In 1738 he met and befriended then Crown
Prince Fredrick of Prussia. Like Fredrick, Bielfeld was a freemason. With
Fredrick's ascent to the Crown in 1740, Bielfeld started his diplomatic
career as Counselor to the Prussian Consulates in Hannover, later in 
London and Berlin. In 1745 he became tutor to Prince August Ferdinand, in
1747 he became curator of Prussia's universities and director of Berlin's
famous Charit\'e hospital. In 1748 Bielfeld was ennobled as Baron, and
after 15 years of service to Prussia --at the age of 38-- he withdrew
from a \emph{vita attiva} to a \emph{vita contemplativa} 
at his properties in Altenburg in
the Eastern part of Germany in 1755. Due to the Seven Years' War 
Bielfeld had to leave Altenburg in 1757 for his native Hamburg, returning only
in 1763. He died in Altenburg in 1770, at the age of 53
[Reinert 2013:7].
\end{quotation}
The primary source for Bielfeld's biography is the eulogy,
read on May 30, 1770 at a public session of the Berlin Academy by the perpetual
secretary Johann Heinrich Samuel Formey (1711-1797). This eulogy was
reproduced as a preface to the posthumous edition of the third volume
of the \emph{Institutions politiques} [Bielfeld 1772:xi-xviii].
In their review of that volume, the authors of the \emph{Biblioth\`eque
des Sciences et des Beaux-Arts} remark:
\begin{quotation}
Cet Eloge est tel que devroient \^etre toutes les
Pi\`eces de ce genre. L'habile \& judicieux Secretaire
y rend une exacte justice \`a l'Acad\'emicien dont
il d\'eplore la perte, il c\'el\`ebre son m\'erite
\& ses talens, mais sans flatterie, sans enthousiasme; \& en
caract\'erisant, en appr\'eciant ses Ouvrages, il ne 
dissimule pas ce qu'ils peuvent avoir de
d\'efectueux. [Chais and de Joncourt 1772:133]
\vskip 1mm
[This eulogy is like all works of the kind should be.
The deft and wise secretary, gives exact credit to the academician he
laments the loss of, celebrates his merit and talents, though without
flattery, whithout enthusiasm; and while characterizing, appreciating
his works, he does not conceal their flaws.]
\end{quotation}
Indeed, Formey's review of Bielfeld's work is unusually critical.
\begin{quotation}
C'est ici le lieu de parler de ses productions. Si elles ne sont pas
de la premiere classe, on ne peut leur contester, au moins \`a la
pl\^upart d'entr'elles, un rang plus ou moins honorable dans la
seconde. Il avoit fait son coup d'essai par une Traduction; c'est
celle des \emph{Consid\'erations sur les causes de la grandeur \& de
  la d\'ecadence des Romains}; elle fut imprim\'ee \& bien
accueillie. Il donna ensuite les \emph{Progr\`es des Allemans dans les
  Belles-Lettres}; Ouvrage int\'eressant, mais qui n'est pourtant
qu'esquiss\'e \& assez incorrect. Je serois tent\'e de passer sous
silence ses \emph{Amusemens dramatiques}, qui n'amus\`erent jamais que
lui. Mais ses \emph{Institutions politiques}, comme je l'ai d\'ej\`a
insinu\'e, sont un Livre v\'eritablement estimable. Il n'y est pas
cr\'eateur, mais il n'y est pas non plus simple compilateur. Il a fait
un bon choix; il y a mis un bon ordre; \& ce qui est de lui ne
d\'epare pas ce que des Auteurs distingu\'es peuvent lui avoir
fourni. Un Critique des plus mordans voulut couler \`a fond ce Livre,
mais il n'y r\'eussit pas. Si ses censures \'etoient quelquefois
fond\'ees, leur aigreur g\^atoit tout; M. de Bielfeld, naturellement
doux \& poli, se fit bien plus d'honneur encore par la mod\'eration de
ses r\'eponses que par leur solidit\'e. Ses \emph{Lettres
  famili\`eres} furent un enfant de son loisir: mais un enfant
g\^at\'e \& beaucoup trop familier. Ses \emph{Traits d'\'Erudition
  universelle} ne sont que des traits; l'ensemble y manque: les
jeunes gens peuvent pourtant en tirer quelque parti. Enfin il fit une
feuille p\'eriodique en Allemand, intitul\'ee l'\emph{Hermite}; elle
s'est soutenue trois ans; c'est beaucoup pour ce genre d'Ouvrage, qui
n'a pas la vie longue, pour peu qu'il soit foible. 
[Formey 1770:71]
\vskip 1mm
[This is the place to speak about his productions. If they are not of
the first class, one cannot deny them, or at least most of them, an
honorable place in the second. He had made his trial run by a
translation; that of the \emph{Consid\'erations sur les causes 
de la grandeur \& de la d\'ecadence des Romains}; it was printed and
well greeted. He then gave the \emph{Progr\`es des Allemans dans les
  Belles-Lettres}; an interesting work, though it is only sketched, and
rather incorrect. I should be tempted to keep quiet about his
\emph{Amusemens dramatiques}, that never amused anybody but himself.
But his \emph{Institutions politiques}, as I have already hinted, are
a truly estimable book. He is not a creator there, but neither is he
a mere compiler. He has made a good choice; he has put a good order;
and what is his does not spoil what more distinguished authors may
have supplied to him. A most caustic critic wanted to bring down that
book, but did not succeed. If his criticisms were sometimes
well-founded, their sourness spoilt everything; M. de Bielfeld,
naturally gentle and polite, did himself even more honor by the
moderation of his answers than by their soundess.
His \emph{Lettres famili\`eres} were a child of his leisure 
time: but a spoilt child,
much too familiar. His \emph{Traits d'\'Erudition
  universelle} are nothing but strokes; a bird's-eye view is missing: yet
young people may draw some benefit from it. Finally, he made a
periodical in German, entitled the \emph{Hermit}; it was maintained
for three years; it is quite long for that kind of work, which is not
long-lived, particularly if it is weak.]
\end{quotation}
Thus according to Formey, the \emph{Institutions politiques} were the
only one among Bielfeld's works to be ``truly estimable''. As we shall
see, it enjoyed a considerable diffusion. 
Chapter XIV of volume II contains Bielfeld's view on
political arithmetics, and volume III is a statistik of
Europe. The \emph{Traits d'\'erudition universelle}, were indeed meant as
a textbook for young readers, and Bielfeld never pretended it to be
more than ``traits'' (strokes). 
\begin{quotation}
Ye Studious Youth! do not repay me with ingratitude; do not accuse me of
presumption, nor imagine that I regard this work as a masterpiece of the
human mind, that makes pretentions to immortality. No, my utmost ambition
is to provide you with a useful work. If you shall interleave these sheets
with blank paper; if you shall read them often, and mark down all the 
observations you will make on each subject, during the course of your 
studies, you can scarce possibly avoid acquiring a valuable portion
of erudition. [Bielfeld 1970a]
\end{quotation}
The ``traits'' were translated as ``elements'' in English,
and somewhat exageratedly as ``curso completo'' in Spanish.
Chapter XIII of volume III is entitled ``Statistique'', and
this is the first occurrence of the word in French in 1767, 
then in English through Hooper's translation [Bielfeld 1770b].
\subsection{Bielfeld's critics}
This section focuses on the increasingly negative perception of
Bielfeld and his work, 
from his death until the end of the nineteenth century. 
Critics of the \emph{Amusemens dramatiques} and \emph{Lettres
  famili\`eres} will be examined first, biographies of Bielfeld will
come next. 

Formey's negative opinion of the \emph{Amusements Dramatiques} (that
``never amused anybody but himself'')
was apparently shared at the publication
time.  
Here is an excerpt from a letter sent by  Friedrich Melchior
Grimm (1723-1807) to Denis Diderot (1713-1784), on April 15, 1774.
\begin{quotation}
J'ai lu, il  n'y a pas long-temps, la pr\'eface
que M. le baron de Bielfeld, Allemand, a mise \`a la t\^ete d'un
recueil de ses com\'edies. Apr\`es cette lecture, j'avoue que je
n'ai as eu le courage de lire la moiti\'e d'une sc\`ene d'une de ses
pi\`eces. Il est impossible de parler sur la mati\`ere que nous
venons de traiter, avec plus de d\'eraison que cet auteur n'a fait.
[Grimm and Diderot 1829:133]
\vskip 1mm
[I have read, not long ago, the preface that Baron Bielfeld, German,
put as a header of a collection of his comedies. After that reading, I
confess that I did not have the courage to read one half of a scene of
one of his plays. It is impossible to talk about the matter we have
just treated, with more folly than this author has done.]
\end{quotation}

In 1763, the English anonymous reviewer of the \emph{Lettres
Famili\`eres} does not
hide his disappointment. 
\begin{quotation}
The reputation this noble Author hath acquired by his \emph{Political
  Institutions}, very naturally excited our curiosity with regard to
his epistolary correspondence. After a very fair and candid perusal,
however, of the Letters before us, we must confess ourselves to have
been a little disappointed in the expectations we had formed of them.
[Monthly Review, vol XXVIII, 1763:516]
\end{quotation}
His conclusion shows no mercy:
\begin{quotation}
Matters of greater importance also, oblige us here to 
dismiss these Letters.
[Monthly Review, vol XXVIII, 1763:523]
\end{quotation}
Yet, the \emph{Lettres famili\`eres} had several editions. Their
English translation was used by several historians
as a source for the reign of Frederick the Great. Some, like 
Joseph Towers (1737-1799) or John Abbott (1805-1877), 
accept Bielfeld's Letters
as a primary source coming from a direct witness, with no negative
judgement [Towers 1788; Abbott 1871].
Others, like Thomas Carlyle (1795-1881), extend their appraisal to
Bielfeld's personality. 
\begin{quotation}
Fantastic Bielfeld, in his semi-fabulous style, has a 
LETTER on this bombardment, attractive to Lovers of the 
Picturesque, (written long afterwards, and dated \&c. WRONG). 
As Bielfeld is a rapid clever creature of the coxcomb sort, 
and doubtless did see Neisse Siege, and entertained seemingly 
a blazing incorrect recollection of it, his Pseudo-Neisse 
Letter may be worth giving, to represent approximately what 
kind of scene it was there at Neisse in the October nights.
[Carlyle 1865]
\end{quotation}
Fifteen years later, Andrew Hamilton, 
in his history of the Rheinsberg court devotes chapter XVI of volume I to
Bielfeld's letters. Not only
does he dismiss them as a reliable source, 
but he goes even farther in caricature than Carlyle.
\begin{quotation}
The man himself was certainly not a very nice sort of man, and the
conditions of is life brought into play usually not the best of what
was in him, but oftentimes rather the worst. He was a
\emph{parvenu}, but never quite succeeded in climbing to any of the
heights the reaching of which makes \emph{parvenu]}ship a safe and
honourable calling.

[\ldots]

Bielfeld was not found full grown in any post. About 1755, I
think, after some \emph{d\'em\^el\'es} with the King, he left
Prussia altogether, and settled -- having got married to a rich wife
in the meanwhile -- on an estate in Altenburg. There, and in
Hamburgh, he spent the rest of his life, writing books of no great
value. He was estimable in private life; in the private life of the
province or of his native town all the more estimable no doubt
because of his former greatness, or rather great expectations.
[Hamilton 1880:211-219]
\end{quotation}

We shall now examine some of the biographical sketches that can be
found in 18\textsuperscript{th}- and 19\textsuperscript{th}-century 
dictionaries. All rely on Formey's eulogy,
whether explicitly cited or not. The earliest come from
French dictionaries, that were later translated into other languages.
The first one is the \emph{Dictionnaire Universel des
  Sciences Morale, \'Economique, 
Politique} by Jean-Baptiste Robinet (1735-1820). There, 
Bielfeld has a laudatory forty-five pages entry
where Formey's eulogy is almost completely copied out (except
the criticisms on Bielfeld's works). The three volumes of the 
\emph{Institutions politiques}, are extensively quoted [Robinet
1779:248-293]. The next dictionary is
the \emph{Biographie Universelle} by 
François-Xavier de Feller (1735-1802), which had many successive editions
from 1781 on, and was continued until 1850, long after de Feller's
death. Bielfeld is absent from the first edition in
1781; he is present in all successive editions from 1784 until 1850.
The sketch is very critical, and strongly
biased by de Feller's viewpoint as a catholic clergyman.
For instance about [Bielfeld 1772]:
\begin{quotation}
On y trouve une description g\'eographique de l'Europe, m\^el\'ee de
r\'eflexions politiques; il est facile de voir, en lisant les articles
qui concernent l'Espagne, le Portugal, l'Italie, etc., qu'il \'ecrit
en bon protestant. On y lit des choses d'une fausset\'e \'evidente,
que la passion seule lui a dict\'ees. [de Feller 1784:708]
\vskip 1mm
[There can be found a geographic description of Europe, mixed with
political reflections; it is easy to see, by reading the articles
concerning Spain, Portugal, Italy, etc.; that he writes as a good
Protestant. Obviously false things can be read there, that only passion
has dictated to him.]
\end{quotation}
About \emph{Progr\`es des Allemans dans les belles-lettres},
de Feller writes:
\begin{quotation}
Mauvaise compilation o\`u le fanatisme
  protestant tient souvent lieu de critique. Si on devait juger des
  progr\`es de la civilisation et des sciences chez les Allemands, par
  la mani\`ere dont son livre est r\'edig\'e, il n'y aurait point de
  nation en Europe moins avanc\'ee. [de Feller 1784:708]
\vskip 1mm
[Bad compilation where Protestant fanatism often serves as a criticism. If
progress of civilization and science among the German had to be judged
by the way his book is written, there would not be any less advanced
nation in Europe.]
\end{quotation}
Rather hypocritically, in his first version,
de Feller hides behind Formey's eulogy: 
\begin{quotation}
Ce que nous disons des ouvrages de Bielfeld est presque tir\'e mot \`a
mot de son \'eloge, fait par un de ses intimes amis, \& l\^u dans une
assembl\'ee publique de l'acad\'emie de Berlin, en 1770. [de Feller
1784:708].
\vskip 1mm
[What we say about Bielfeld's works comes almost verbatim from his
eulogy, made by one of his close friends, and read in a public
assembly of the Academy of Berlin, in 1770.]
\end{quotation}
From 1795 on, the conclusion of the article has changed:
\begin{quotation}
Un de ses intimes amis a lu son ``Eloge'' dans une assembl\'ee
publique de l'acad\'emie de Berlin, en 1770: 
on comprend bien que l'auteur et ses
ouvrages n'y sont pas s\'ev\`erement jug\'es. [de Feller 1797:225]
\vskip 1mm
[One of his close friends has read his ``Eulogy'' in a public assembly
of the academy of Berlin, in 1770: it is easily understood that the
author and his works are not severely judged.] 
\end{quotation}
That version was translated into Spanish [Oliva 1830:537].
Fran\c{c}ois Guizot (1787-1874) and Adrien Jean Quentin Beuchot
(1773-1851) signed Bielfeld's sketch in Michaud's \emph{Biographie Universelle}
[1811]. Chalmers' biography is a
translation, from which the following is quoted.
It was also translated into Italian (page 122 of 
\emph{Biografia universale antica e moderna}, volume VI, Venezia:
Missaglia, 1822).
\begin{quotation}
In a journey which he made to Brunswick, he became acquainted
with Frederic II. then prince royal, who, on coming to the throne, took
him into his service, and sent him, as secretary of legation, with
count de Truchses, Prussian ambassador to the court of St. Jame's, but
discovering that the baron's talents were not calculated for
diplomatic affairs, he, in 1745, appointed him preceptor to prince
Augustus Ferdinand his brother;

[\ldots]

He wrote 1. ``Institutions politiques,'' [\ldots], 
the only work  from his pen that
retained its reputation on the continent. 
Even the empress Catherine
II. of Russia, condescended to write notes on it.
[Chalmers 1812:250]
\end{quotation}
Until 1830, substantial sketches of Bielfeld life and works 
can be found in most dictionaries of biography;
after 1830 the entries become increasingly short, neutral, and
sometimes erroneous. Bielfeld is presented alternatively as 
a `celebrated modern writer' [Maunder 1838:92], a `publicist'
[Hoefer 1853:29], or
a `German statesman' [Hyamson 1916:60]. His birth date is often
cited as 1716, or even `around 1712' [Thomas 1870:355].


%
\section{Bielfeld and political arithmetic}
%
\subsection{Of political calculations}
\begin{figure}[!ht]
\centerline{
\includegraphics[width=12cm]{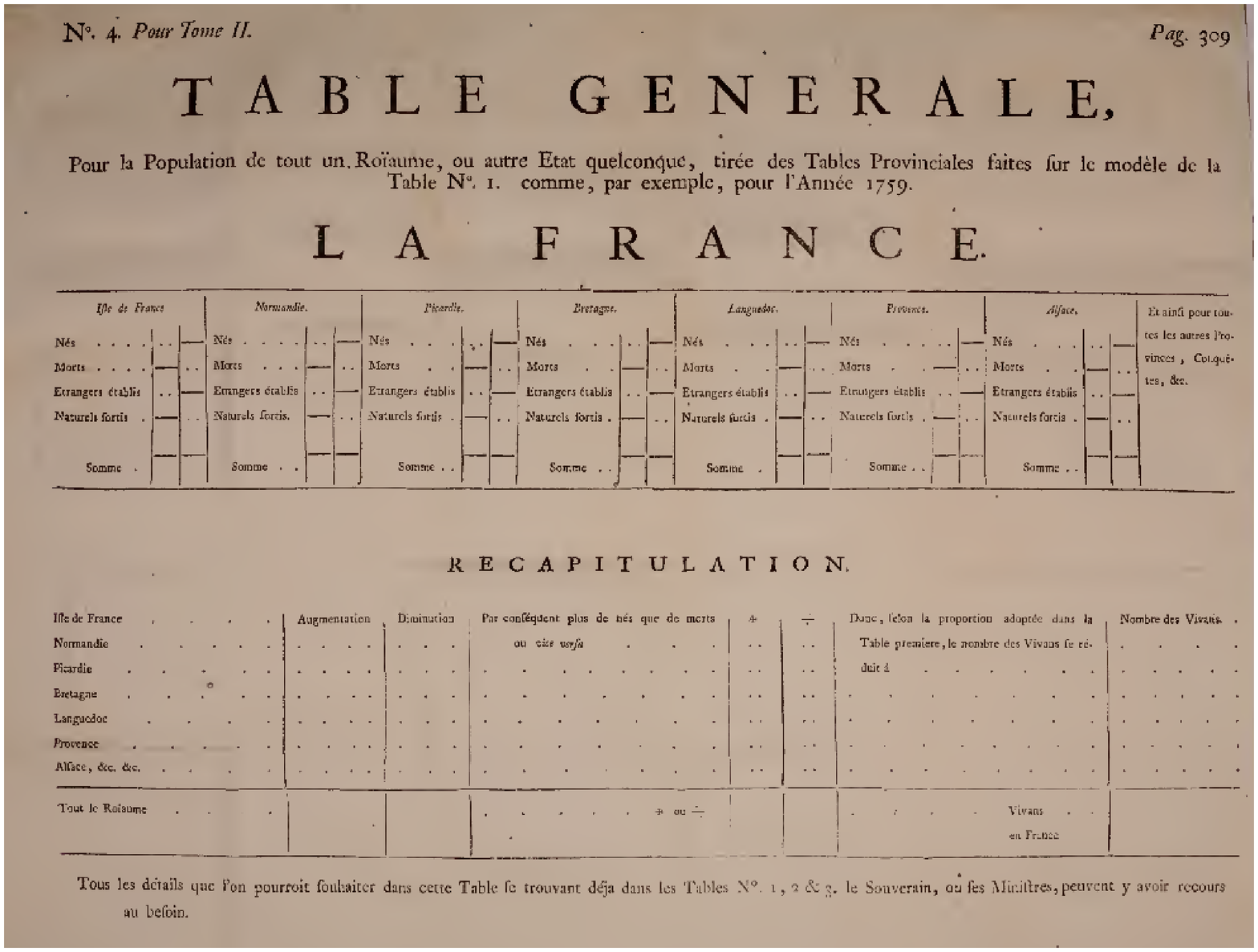}
} 
\caption{Model of recap chart for the population of a country
  [Bielfeld 1760b:309].}
\label{fig:Recap}
\end{figure}

This section reviews Chapter XIV, Volume II of Bielfeld's \emph{Institutions
  politiques}, entitled ``Des Calculs Politiques'' [Bielfeld
1760b:263-309]. As Bielfeld precised in the introduction to Volume I,
the \emph{Institutions politiques} have 
no pretention to innovation. 
\begin{quotation}
[\ldots] on ne doit
pas s'attendre \`a trouver dans tout le cours de cet Ouvrage des
id\'ees nouvelles que personne n'a eues, des d\'ecouvertes
singuli\`eres qui sont le fruit d'une imagination brillante.
[Bielfeld 1760a:6]
\vskip 1mm
[One cannot expect to find in all this work new ideas that no one has
had, singular discoveries springing from a brilliant imagination.]
\end{quotation}
The only data presented essentially come from the mortality tables of
Wilhem Kerseboom (1691-1771) [1738-42].
In the first lines of Chapter XIV, Bielfeld clearly states his goals.
\begin{quotation}
L'Arithm\'etique Politique a \'et\'e r\'eduite, depuis environ soi\-xante \& dix
ans, en Science particuli\`ere. Des Calculateurs habiles \& infatigables
se sont fortement appliqu\'es \`a la perfectionner; \& leurs Ouvrages ont
contribu\'e \`a la rendre si c\'el\`ebre, qu'aujourd'hui les grands Hommes
d'Etat semblent \^etre dans l'opinion qu'elle est indispensablement
n\'ecessaire pour r\'egir un pa\"\i s. Ces consid\'erations nous obligent d'en
faire quelque mention dans cet Ouvrage. Nous t\^acherons d'indiquer (1)
l'origine et l'histoire de cette Science (2) les objets sur lesquels
elle peut porter (3) le degr\'e de certitude dont elle est susceptible,
(4) \`a quel point elle est applicable dans la pratique du Gouvernement,
(5) les Principes sur lesquels elle se fonde \& (6) les op\'erations
qu'elle emplo\"\i e pour d\'ecouvrir ce qu'elle cherche: Car entrer dans
les Calculs m\^emes, ou les appliquer \`a divers pa\"\i s de l'Europe, ce
serait s'engager dans une entreprise trop vaste pour les bornes de
notre Plan, \& nous ne pourrions que transcrire ce que tout Lecteur
peut trouver dans les Auteurs qui ont trait\'e Sist\'ematiquement cette
matiere, \& que nous citerons chemin faisant.
[Bielfeld 1760b:263]
\vskip 1mm
[Political arithmetic has been reduced for about seventy years, into
a particular science. Deft and tireless calculators have devoted
themselves to improving it; and their works have contributed to making
it so famous, that nowadays great statesmen seem to 
believe that it is indispensably necessary to rule a country.
These considerations oblige us to mention it in this publication.
We shall try to indicate (1) the origin and history of this science
(2) the objects on which it can bear (3) the degree of certainty it
can attain (4) up to which point it is applicable in the practice of
government (5) the principles upon which it is grounded (6) the
operations it uses to find what it looks for: because entering into the
very calculations, or applying them to various countries of Europe,
would be embarking into too vast an endeavor for the limits of our
plan, and we could only transcribe what any reader can find in the
authors that have systematically treated the matter, and that we shall
cite along the way.]
\end{quotation}
As announced, the chapter begins with a historical review of the
field, introduced by an emphatic tribute to the English forerunners. 
\begin{quotation}
L'Arithm\'etique Politique est n\'ee dans le terroir qui devoit
naturellement la produire, c'est-\`a-dire en Angleterre. Un pa\"\i s ou
toutes les parties des Math\'ematiques sont cultiv\'ees avec tant de soin,
qui a l'honneur de l'invention de tant de Calculs fameux, qui a
produit le c\'el\`ebre Newton, Pere de tous les Calculs; une Nation qui
s\c{c}ait peser jusqu'aux Astres, \& qui joint \`a ce talent un go\^ut d\'ecid\'e
pour la Politique, ne pouvoit manquer de r\'eduire aux Principes du
Calcul les objets principaux du Gouvernement de l'Etat.
[Bielfeld 1760b:263]
\vskip 1mm
[Political arithmetic was born in the land that was naturally meant to
produce it, \emph{i.e.} in England. A country where all parts of
Mathematics are cultivated with so much care; that has the honor of
the invention of so many famous calculations, that has produced the
celebrated Newton, father of all calculations; a nation who can weight
up to the stars, and which joins to that talent a decided taste for
politic, could not fail to reduce to principles of calculation, the
main objects of the government of a state.]
\end{quotation}
The list of authors from all over Europe cited by Bielfeld,
evidences his broad knowledge of the field: see Th\'er\'e and
Rohrbasser in [Martin 2003:304-7]. With Bielfeld's spelling, his
sources come from
England (Graunt, Petty, Derham, de Moivre, Halley, King, Arbuthnot,
Hogdson, Maitland, Simpson, Hume), 
France (Vauban, de St Pierre, du Tot, Melon, Davenant,
Desparcieux, Buffon), 
Holland (Nieuwenhyt, Struyk, Kersseboom, 's Gravesande), 
Germany (Susmilch, Justi, Kundmann),
Switzerland (Bernouilli), Spain (de Uztaritz),
Sweden (Fayot, Wargulin, Berch).
Thorough accounts of the most important works, such has those of
Petty, S\"ussmilch, Vauban, Kersseboom, etc. are given, and the text
contains some
personal appreciations:   
\begin{quotation}
En Suisse, l'infatigable \emph{M. Bernouilli} a \'eclairci divers points
relatifs au Calcul Politique, \& l'on s\c{c}ait quel est le juste cas que
l'on doit faire de tout ce qui sort de la plume de ce grand homme en
qui la Science des Math\'ematiques semble \^etre inn\'ee.
[Bielfeld 1760b:273]
\vskip 1mm
[In Switzerland, the tireless Mr. Bernoulli has cleared up
different points relative to political calculation, and one knows 
how important a case must be made of what comes from the pen of that great
man in which the science of mathematics seems to be innate.]
\end{quotation}
After listing the type of data political arithmetic deals with (census,
taxes, income, etc.),
Bielfeld addresses the crucial issue of precision.
\begin{quotation}
On serait trop heureux si tous les diff\'erens Calculs Politiques,
dont nous avons parl\'e jusqu'ici, pouvoient se faire avec une
pr\'ecision parfaite; mais il s'en faut de beaucoup qu'ils soient
susceptibles d'une certitude math\'ematique.
[Bielfeld 1760b:282]
\vskip 1mm
[One would be too fortunate if all political calculations, that have
been spoken of so far, could be done with a perfect precision; but
they are far from being subject to mathematical certainty.]
\end{quotation}
His argument is twofold. On the one hand, errors in 
 counting are inevitable.  
\begin{quotation}
La m\^eme incertitude r\`egne dans les D\'enombremens.
Quelque soin qu'on prenne, il est impossible de les faire
avec une enti\`ere pr\'ecision. On ne compte pas les hommes,
non plus que les feuilles d'une Foret, ou que tous les
\^etres qui se changent \& se renouvellent sans cesse.
Chaque Ville ressemble en cela \`a un Colombier, ou \`a une
Ruche ouverte, dont les habitans to\^ujours en mouvement entrent,
sortent, s'agitent sans rel\^ache, \& confondroient l'exactitude du
Calculateur le plus infatigable qui voudroit d\'eterminer leur 
nombre.
[Bielfeld 1760b:287]
\vskip 1mm
[The same uncertainty reigns in counting. Whatever care
is taken, it is impossible to make them with a thorough precision. One
cannot count men, no more than leaves in a forest, or all beings that
change and renew themselves incessantly. 
Each city resembles in that respect
a dovecot, or an open beehive, the inhabitants of which, always
moving, go in and out, move relentlessly, and would confuse the
accuracy of the most tireless calculator who would want to determine
their number.]  
\end{quotation}
On the other hand, absolute precision is a useless mathematical
fantasy (see [Bru 1988:11; Feldman 2005:14]). 
\begin{quotation}
Mais la politique n'a pas besoin, dans cette affaire-ci, d'une
pareille certitude. Elle peut se contenter tr\`es bien d'une th\'eorie
vraisemblable sur tous ces objets, pourvu que cette th\'eorie soit aussi
approchante de la v\'erit\'e qu'il est possible; \& c'est ce \`a quoi
tendent tous les efforts des Calculateurs Politiques, qui devroient
\^etre second\'es, dans les pays bien polic\'es, par le Gouvernement
m\^eme. Vouloir aller au-del\`a, \& pr\'etendre \`a la 
pr\'ecision Math\'ematique
dans cette mati\`ere, ce seroit chercher un objet de sp\'eculation et de
curiosit\'e, comme la Quadrature du Cercle.
[Bielfeld 1760b:288]
\vskip 1mm
[However, politic does not need, in this matter, such a precision. It
can make do with a likely theory on all these objects, provided that
theory is as close to truth as possible; and this is what all efforts
of political calculators tend to, and they should be assisted in all
well civilized countries, by the government itself. Wanting to go
beyond, and yearning for mathematical precision in that matter,
would be like searching an object of speculation and curiosity, like
squaring the circle.]
\end{quotation}
Furthermore, Bielfeld repeatedly warns the reader against too strict
an application of mathematics to human affairs, sometimes with a touch
of humour.
\begin{quotation}
Un Roi, par exemple, dans son Conseil, un Ministre,
dans son Cabinet, qui calculeroit, comme on l'a dit,
les affaires par $a+b\div c$
courroit risque de prendre \`a tout moment une r\'esolution
qui seroit \'egale \`a z\'ero.
[Bielfeld 1760b:288]
\vskip 1mm
[A king for instance, in his council, a minister, in his cabinet, who
would calculate, as has been said, the affairs by $a+b\div c$ would be
at risk to take at any time a resolution that would be equal to zero.]
\end{quotation}
Even though absolute precision is illusory, and not even necessary,
Bielfeld still advocates accurate data collection.
\begin{quotation}
Il faut le rep\'eter encore ici, les Listes des Enfans batis\'es, des
Mariages, \& des Morts forment la baze de toute cette
Arithm\'etique. Il est donc n\'ecessaire que le Souverain ordonne \`a
tous les Cur\'es, tant des Villes que de la Campagne, sans exception,
de tenir des Registres exacts de tout ce qui se passe relativement \`a
cette mati\`ere dans l'\'etendu\"e de leurs Paroisses respectives. Une
seule Paroisse qui manque rend tout le Calcul imparfait \& faux.
[Bielfeld 1760b:289]
\vskip 1mm
[This must be repeated here again, the lists of baptised children, of
marriages, and of deaths, form the basis of all this arithmetic. Thus it
is necessary that the sovereign command all ministers, in cities as
well as in the countryside, to keep exact registers of what happens in
these matters in the whole extent of their respective parishes. A single
parish missing makes all the calculation imperfect and false.] 
\end{quotation}
After the data have been carefully collected, information must be
extracted.
\begin{quotation}
Lorsque le Souverain s'est procur\'e toutes ces Dates, le plus
exactement qu'il est possible, il peut donner de l'occupation \`a ses
Calculateurs Politiques, qui en font des r\'esum\'es tr\`es
instructifs en calculant\ldots
[Bielfeld 1760b:290]
\vskip 1mm
[When the sovereign has obtained all these data, as exactly as
possible, he can give occupation to his political calculators who will
make very instructive summaries by calculating\ldots] 
\end{quotation}
Bielfeld is perfectly aware that data collection has a long tradition,
yet he is not completely satisfied with the method.
\begin{quotation}
Nous en trouvons des mod\`eles dans celles qui ont \'et\'e faites \`a
Londres \& \`a Vienne, \& qui sont rapport\'ees par
M\textsuperscript{rs.} Graunt, Kundmann, Susmilch, \&c. Mais ces
mod\`eles souffrent encore, ce me semble, des rectifications \& des
augmentations que les ordres du Souverain, ou de ses Ministres,
peuvent leur donner tr\`es facilement. Il me paroit n\'ecessaire de
joindre \`a ce Paragraphe (Voyez \`a la fin du Chapitre), les
mod\`eles de quatre Tables qui renferment les parties de la
Population, ou les Dates, les plus essentielles pour l'Arithm\'etique
Politique en g\'en\'eral, \& dont la seule inspection pourra fournir
une id\'ee un peu plus claire de ce genre de Calcul.
[Bielfeld 1760b:291]
\vskip 1mm
[Models can be found in those that have been done in London and Vienna,
and that have been reported by Mr Graunt, Kundmann, Susmilch, etc. But
these models are still missing, I think, some rectifications and expansions
that orders of the sovereign or his ministers may easily give them.
It seems necessary to me to
append to this paragraph (see at the end of the chapter) the models of
four tables which enclose the parts of the populattion, or the data,
which are most essential for political arithmetic in general, and of
which the mere inspection will give a clearer idea of that kind of
calculation.]
\end{quotation}
Indeed, four models of tables are inserted at the end of the
Chapter. The first three are at the province scale. The first one
gives population counts
according to categories, 
the second one is a mortality table (deaths per age group), the third
one is a casualty table (deaths per cause), following the English model. 
The fourth and last table (Figure \ref{fig:Recap}) is a recap
chart for the whole country. It is interesting to notice that such
models of tables, and their generalization to all sorts of data, 
were the basic tool of the ``great avalanche of printed
numbers'', half a century later [Hacking 1990:2]. Quite similar charts
were submitted by Sinclair to the Scottish ministers
and later helped him summarizing the wealth of data they 
had sent [Sinclair 1791:viii; 1798:xv ff.].

Thus, with the help of his tables, specialized calculators could
provide the governement with all sorts of usefull quantified
information: 
Bielfeld anticipates on the ``bureaus'' that would later
flourish all over Europe [Hacking 1990:27].
\begin{quotation}
Il me semble que la chose vaudroit bien la peine que chaque Etat
entretint un ou deux Calculateurs Politiques qui n'eussent qu'\`a
s'appliquer uniquement \`a recueillir les Listes des Paroisses, \& \`a en
former des Tables telles que je viens de les proposer. Cette l\'eg\`ere
d\'epense seroit bient\^ot compens\'ee par les avantages que l'Etat
retireroit de la justesse de toutes les mesures qu'on pourroit prendre
dans le Gouvernement.
[Bielfeld 1760b:291]
\vskip 1mm
[It seems to me that it would be worth the trouble 
for each state to maintain one
or two political calculators whose only ocupation would be to collect
the lists from parishes, and build tables such as those proposed
above. This slight expenditure, would be compensated soon by the
advantages that the state would get out from the accuracy of all measures
that could be taken in the government.] 
\end{quotation} 
Thus the methodology of data collection that was to be generalized at 
the turn of the century, is clearly described. What about the theory?
The analogy with games of chance is blatant to Bielfeld, as it must have
been to his readers. About life insurance, he says:
\begin{quotation}
Il est clair que cet \'etablissement, ou ce Contrat, n'est dans le fond
qu'un \emph{Jeu de Hazard}, qui a ses \emph{chances} comme tous les autres,
\& dont le profit, ou la perte, roule uniquement sur la dur\'ee, plus ou
moins longue de la vie de de celui sur la tete duquel le fond perdu est plac\'e.
[Bielfeld 1760b:297]
\vskip 1mm
[It is clear that this establishment, or this contract, is but a game,
which has its chances as any other, and the profit or loss of which, depends 
only on the duration, more or less long, of the life of the person on which
the unredeemable fund is placed.]
\end{quotation}
The obvious corollary is equity of expectated profits.
\begin{quotation}
Nous venons de remarquer que tout \'etablissement de cette nature
a ses chances; il doit donc avoir aussi ses justes \emph{proportions},
si une des partie contractantes, ou, si vous voulez, des Joueurs, ne doit
pas \^etre une dupe manifeste.
[Bielfeld 1760b:301]
\vskip 1mm
[We have just remarked that any establisment of this nature has
its chances; therefore it must have also its fair proportions,
if one of the contracting parties, or of the gamblers if you want,
is not to be a patent dupe.]
\end{quotation}
Bielfeld has read Bernoulli, and expresses his deepest admiration.
Yet he does not make clear that he has fully grasped the law of large numbers, 
in which he professes an almost candid belief.
\begin{quotation}
La fixation de toutes les proportions des diff\'erens \^ages de la vie
humaine sert aussi de r\`egle pour le Calcul des \emph{Tontines}: Car,
je l'ai dit plus d'une fois, il y a dans la Mortalit\'e d'un grand
nombre d'hommes, \& dans une r\'evolution consid\'erable d'ann\'ees, moins de
hazard, moins d'irr\'egularit\'e qu'on ne pense. Les m\^emes chances, les
m\^emes proportions revienent toujours au bout de quelque tems. 
Il semble que tout ce qui est dans la Nature soit asservi \`a de certaines
r\'evolutions p\'eriodiques \& presque uniformes dans un grand espace de tems.
Dans les Jeux de hazard m\^eme, ce qu'on appelle proprement hazard, ou fortune, y
entre pour beaucoup moins qu'on ne s'imagine, ou du moins cette
fortune n'y est pas si bizarre, si constante, ni si inconstante,
qu'elle paroit l'\^etre au premier coup-d'\oe il, ainsi que l'a tr\`es bien
d\'emontr\'e le c\'el\`ebre Mr. Bernouilli, \& quelques autres
Calculateurs. Si l'on d\'ecompte de ces coups du hazard, les tours
d'adresse, les tromperies, \& les fautes de jugement qui se font \`a de
pareils Jeux, on  verra qu'au bout de quelques ann\'ees des Joueurs qui
luttent souvent l'un contre l'autre balanceront assez leur gain ou
leur perte, ce qui prouve clairement l'\'equilibre dans les retours de
la fortune.
[Bielfeld 1760b:305]
\vskip 1mm
[Fixing the proportions of the different ages of human life
serves also as a rule for the calculation of tontines. Because,
I have said it more than once, there is in the mortality of a large
number of men, less hazard, less irregularity than one would think.
The same chances, the same proportions, always come back after some time.
It seems that all that exists in nature is subject to certain
periodic revolutions, almost uniform in a large interval of time.
In hazard games themselves, what is properly called hazard, or 
fortune, enters for much less as one would imagine, or at least that
fortune is not as bizarre, nor as constant, nor as inconstant, that
it seems to be at first sight, as have so well shown the famous Mr. 
Bernouilli and some other calculators. If from those hazard strokes, 
the skilful tricks, frauds, or judgement mistakes that occur in 
such games are discounted, one will see that after some years,
players that often gamble against each other, will rather 
balance their profit or loss, which clearly proves equilibrium in 
the fluctuations of fortune.]
\end{quotation}
Bielfeld's strong faith in the law of large numbers seems to
announce the ``Constants of
Nature and Arts'', dear to Babbage in the following century
[Hacking 1990:52]. It is probably the cause of his erroneous conjecture
about the stability of human population.
\begin{quotation}
Enfin on serait tent\'e de croire que la quantit\'e d'hommes r\'epandus
sur la surface de la Terre a presque toujours \'et\'e la m\^eme,
ainsi que celle des autres Cr\'eatures, \& que ces choses pourroient
bien se soutenir dans le m\^eme arrangement jusqu'\`a la fin des Si\`ecles.
[Bielfeld 1760b:285]
\vskip 1mm
[At last, one would be tempted to believe that the number of men across
the Earth, has almost always been the same, as weel as that of other 
creatures, and that these things could well last in the same arrangement 
until the end of centuries.] 
\end{quotation}
Spengler [1942:chapter III], has remarked
that Bielfeld's belief in population stability somehow contradicts his
repeated recommendation to take any measure susceptible to
increase the population of a country, which he carefully lists 
[Bielfeld 1760b:294-296].  
As it turned out, that mistake
focused many of S\"ussmilch's attacks 
(S\"ussmilch is the `most caustic critic' alluded to by Formey): see
[S\"ussmilch 1741-77:129,377]

Overall, it can be observed that Bielfeld's views
of political arithmetics appears
quite coherent with
other contemporaneous accounts, in particular those of Condorcet in the 
\emph{Encyclop\'edie M\'ethodique} [Bru and Cr\'epel 1994, Feldmann 2005].
\subsection{Diffusion in Europe}
This section firstly examines the diffusion of Bielfeld's
\emph{Institutions politiques}, then focuses on the impact of the
chapter on political arithmetics. 
The controversy with S\"ussmilch and its indirect
effect on the history of the discipline through the works of Euler is
treated last.

Among the three volumes [Bielfeld 1760a; 1760b; 1772], the first two
had an immediate impact, much larger than the posthumous third. 
They were the subject of two extensive and
elogious reviews in the \emph{Biblioth\`eque des Sciences et des Arts}
by Charles Chais and \'Elie de Joncourt [1760a:66-90; 1760b:267-288].
An interesting trace of personal diffusion can be found in the
correspondence of Francesco Algarotti (1712-1764) with Voltaire.
\begin{quotation}
Ho ricevuto jeri una lettera del nostro amico Formey;
nella quale egli mi dice: \emph{Que dites-vous de Bielfeld chevalier 
de s. Anne pour avoir fait des} institutions politiques,
\emph{qui effacent Montesquieu~?} 
[Algarotti 1764:172]
\vskip 1mm
[I have received yesterday a letter from our friend Formey, in which
he tells me: ``what do say of Bielfeld, knigth of St. Ann, for having
made political institutions that upstage Montesquieu?]
\end{quotation}
The compliment is somewhat exaggerated, and it might be 
tainted with some trace of irony, the correspondent being Voltaire;
yet as we have seen in his Eulogy, Formey had a sincere appreciation of
the book. Moreover, the allusion to Bielfeld's recent elevation to
the order of St. Ann is a way of recalling the praise of the
Empress Catherine II. of Russia. Some have said that 
she had ``condescended to write notes on it''
[Chalmers 1812:250], others that she had ``placed them, covered with
her own handwritten notes, 
in her library, next to the \emph{Esprit des lois}''
[Sayous 1861:361]. We now know that Catherine II. not only annotated
the political institutions, she 
\begin{quotation}
consulted them extensively in preparing her \emph{great
  Instruction}, especially in Chapter 21 on police which was only
published in February 1768, the same year that the first volume of an
official Russian translation of Bielfeld's treatise appeared.
[Alexander 1988:104]
\end{quotation}
Evidence for the success of the \emph{Institutions politiques} also
come from debates. Some pamphlets were written to refute one single
point of Bielfeld: for instance Amor\'os [1777] on the paternal
consent for marriage, or de J\'ocano y Madaria [1793] on the
double-entry accounting system. Controversies, fueled by religious
passions, sometimes raged accross borders
[Mu\~{n}oz 1778]. Bielfeld is extensively quoted from Portugal
[da Cunha 1794] to Italy [Gautieri 1804]. 
In his treatise of 1823, Zambelli, an
italian translator of Bentham, often cites Bielfeld of which he says:
\begin{quotation}
Infatti il barone di \emph{Bielfeld} nella sua grand'opera
delle \emph{Instituzioni politiche} ha abbraciato tuttoci\`o, che in
latissimo sense pu\`o entrare nell'amministrazione  di uno stato.
[Zambelli 1823:317]
\vskip 1mm
[Indeed, the Baron Bielfeld, in his great work \emph{political
  institutions}, has embraced everything that, in the widest sense,
can enter the administration of a state.]
\end{quotation}
Thus Carpenter's selection of the \emph{Institutions politiques} as one of
the forty `economic bestsellers before 1850' [Carpenter
1975:17] is fully justified. Why then did they sink into oblivion
after 1850? Pierre-Andr\'e Sayous has a sensible explanation:
\begin{quotation}
L'ouvrage qui est \'eclos sous cette influence, a familiaris\'e les
souverains du Nord avec les id\'ees de justice et d'administration
bienfaisantes pour les peuples, et aujourd'hui encore que toutes ces
id\'ees, devenues des lieux communs ou plut\^ot des principes consacr\'es et
pratiqu\'es, ont rendu le livre inutile, les \emph{Institutions
politiques} du baron de Bielfeld conservent encore un grand m\'erite.
[Sayous 1861:363]
\vskip 1mm
[The work that was born under that influence, has familiarized the
sovereigns of the North with the ideas of justice and administration
beneficial to the peoples, and nowadays when all those ideas, having
become commonplace or rather established and practical principles,
have rended the book useless, the  \emph{political institutions} of
Baron Bielfeld, still retain a great merit.]
\end{quotation}
What about Chapter XIV?
As remarked by Th\'er\'e and Rohrbasser in [Martin 2003:304-307],
Bielfeld is not cited in the rather comprehensive
Italian review of political arithmetic by
Gaeta and Fontana [Cr\'epel 2003b]. Does this mean that his writings
went completely unnoticed? Evidence of the contrary can be found in
several editions of popular encyclopedias, mainly those of Robinet and
Felice [Cr\'epel 2003a]. We have seen that the first one has an
extensive biographical 
entry for Bielfeld. Its article \emph{Arithm\'etique Politique}
[Robinet 1778:126-176] is copied from [Bielfeld 1760b]. So are the
articles \emph{Tontine} [Robinet 1783:180-182] and
\emph{Viag\`ere (Rente)} [Robinet 1783:608-616]. Plagiarism was in the
manners of the time, and the fact that Bielfeld [1760b] is explicitly cited
at the end of the article \emph{Tontine} [Cr\'epel 2003a:63] is rather
unusual.  Fortunato Bartolommeo Felice (1723-1789) did not have the
same scruple: he merely signed
Bielfeld's writings with his own initials. The entry
\emph{Arithm\'etique Politique} of his \emph{Code de l'Humanit\'e}
[Felice 1778:512-538] is an exact copy of [Bielfeld 1760b:263-297].
Cr\'epel [2003a:62] cites Diderot about that article, but Hecht [1987:76] 
correctly attributes it to Bielfeld. 

Bielfeld's chapter on political arithmetics was one of many
contemporaneous texts on the subject, with no particular originality,
except perhaps a clearness that might explain why it was preferred to
others by Robinet and Felice. It was certainly
read, remarked, and circulated. One of the few histories of probability which
cites Bielfeld, mentions it as follows: 
\begin{quotation}
Ils furent tous accueillis cependant avec une grande curiosit\'e,
et les principaux d'entre eux, ceux par exemple de Short et de
Corbyn-Morris en Angleterre, de Kerseboom en Hollande, de Sussmilch et
de Bielfeld en Prusse, et chez nous de Deparcieux et de Dupr\'e de
Saint-Maur, m\'erit\`erent v\'eritablement cette faveur.
[Gouraud 1848:42]
\vskip 1mm
[They were all greeted with great curiosity, and the most important of
them, for instance those of Short and Corbyn-Morris in England,
Kersebooom in Holland, Sussmilch and Bielfeld in Prussia, and at home
of Deparcieux and Dupr\'e de Saint-maur, did indeed deserve that favor.]
\end{quotation}
S\"ussmilch would certainly not have liked being paired with
Bielfeld as a political arithmetic author. 
Bielfeld's bold assertions on the stability of human population, added
to some doubts he had cast on the accuracy of S\"ussmilch data
collection, had enfuriated the later: see Hecht's account in
[S\"ussmilch 1741-77:129-130]. The controversy had an indirect, but
long-lasting effect on the history of science. Probably not confident
enough with his grasp of mathematics, S\"ussmilch sought for the help
of Euler, who was generally regarded as the best
mathematician of the time. Bielfeld, who had sitted with Euler at the
foundation of the Academy of Berlin [Laudin 2009:31], would certainly
not contest Euler's authority:
\begin{quotation}  
M. Euler, Math\'ematicien du premier ordre, \& peut-\^etre le
plus grand calculateur qu'il y eut jamais, quitta St. Petersbourg
pour venir s'\'etablir \`a Berlin.
[Bielfeld 1763:336]
\vskip 1mm 
[Mr. Euler, a mathematician of the first rank, and perhaps the greatest
algebraist that ever existed, left St. Petersburg to reside in Berlin.]
\end{quotation}
It is generally considered that Euler is the author of the mathematical
parts in the second edition of  S\"ussmilch's \emph{G\"ottliche
  Ordnung} [Smith 1977].  The series of papers that he wrote
in the 1760's about population, life annuities, lotteries, tontines [Todhunter
1865:239-242, Euler 1923:xxiii-xxv], were probably sparked by his
reading of Bielfeld and the discussions with S\"ussmilch. About
mortality tables, the data that Euler uses are the same as in
[Bielfeld 1760b]: they
come from Kerseboom. The long term effect of Euler's papers is
analyzed by Klyve as follows:
\begin{quotation}
Regardless of whether Euler explicitly showed all his
calculations in S\"ussmilch's book, the claim that in the long run,
an unchecked population will grow geometrically is of seminal
importance. It was this idea about universal geometric population
growth -- hinted at by Wallace, discovered by Euler, circulated by
S\"ussmilch, applied by Price, exploited by Malthus, and
inspirational to Darwin -- which led directly to Darwin's concept of
natural selection.
[Klyve 2013:20]
\end{quotation}
%
\section{Bielfeld and statistik}
%
\subsection{The elements of universal erudition}
As already mentioned, the third volume of the \emph{Institutions
  politiques} is a statistik of Europe [Bielfeld 1772]. Although
cited as such by Meusel [1790:8] and Heuschling [1845:40], 
it did not have much impact. The
probable explanation is that it lacked originality. The information
essentially came from previous similar books by Anton Friedrich B\"usching
(1724-1793), and it was
already outdated at the time of publication. The \emph{Elements of
  universal erudition} had a more long-lasting effect. 
We shall mainly quote from  Volume 3 of 
the English translation [Bielfeld 1770b].

In volume I of the French edition [Bielfeld 1767a], chapter 49, from pages
425 to 518 is devoted to Mathematics, in a very broad acception: it
covers not less than eighteen different sciences which ``extend over
all beings the magnitude of which can be determined by certain
principes, and hence becomes very vast''
(that chapter does not appear in the English translation).
Similarly, chapter XIII of [Bielfeld 1770b:268-279] exposes a very general
definition of statistics, consistent with other German definitions of the
field (from 1784 on, existed
a ``Westphälisches Magazin zur Geographie, Historie und Statistik'';
see also [Meusel 1790:Vorerinnerungen]).
\begin{quotation}
The science, that is called \emph{Statistics, teaches us what is
the political arrangement of all the modern states of the known
world.}
[Bielfeld 1770b:269]
\end{quotation}
At the beginning of the chapter, Bielfeld takes good care to
distinguish statistics from geography, which he treats in Chapter XV.
\begin{quotation}
In geographical treatises, they placed, before the local description
of each country, a sort of account of the principal objects that
composed its system. But these introductions were always imperfect,
naturally very contracted, frequently dubious, and sometimes
absolutely false or ill grounded.
[Bielfeld 1770b:269]
\end{quotation} 
Yet Bielfeld knows what is due to 
B\"usching, even though the ``founder of statistical geography'' 
usually prefers to speak about geography rather
than statistik [B\"usching 1764:Vorbericht].
\begin{quotation}
We must except some of them however, especially those which are to be
found in the excellent geography of M. Busching, an author, whose
assiduity, precision, and discernment, can never be sufficiently
commended.
[Bielfeld 1770b:269]
\end{quotation}
But of course, Bielfeld is aware of Achenwall's founding role.
\begin{quotation}
It would be far from just, in this place, to pass over in silence
the obligations this science has to M. Godfrey Achenwal,
professor at Gottingen, who has not only composed an
Introduction to the political system of the modern 
states of Europe, and another work not less interesting, entitled
Principles of the history of Europe, leading to the knowledge
of the principal states of the present time, but has been also the first
to reduce this important subject into a true system, and has made
a separate science of it, under the title of Statistics, and which he 
professes with great reputation; a science from which history
borrows great lights; which furnishes the best materials for the
constitution of a state, which enriches politics, and which
prepares those of the brightest genius among the studious youth,
to become one day able ministers of the state.
[Bielfeld 1770b:271]
\end{quotation}
The somewhat subtle difference between B\"usching's geography and
Achenwall's statistics is later swept over:
\begin{quotation}
All that occurs in a state is not worthy of remark, but all 
that is worthy of remark in a state, enters necessarily into
statistics.
[Bielfeld 1770b:272]
\end{quotation}
Such a broad definition encompasses
elements of geography, history, economics, politics and, 
almost marginally, political arithmetics.
\begin{quotation}
With regards to the inhabitants, it inquires into their
number and qualities: and for this purpose it makes,
by the aid of political arithmetic, of registers of births and burials,
\&c. the most elaborate and accurate researches possible, into the
number of the inhabitants of a state, and into their genius;
the prevailing character, the industry, the virtues and vices of
a nation.
[Bielfeld 1770b:273]
\end{quotation}
This reminds of the use that B\"usching, even though he does not explicitly
refers to statistics, makes of political arithmetic.
\begin{quotation}
I have set down the probable number of inhabitants in several
countries and great cities, or inserted an account of their births
and burials from the annual Bills of Mortality; but this could not
be done for all.
[B\"usching 1762:vi]
\end{quotation}
After enumerating all the objects of study, Bielfeld addresses the
problem of updating information. He remarks that even
the works ``which approach nearest the exact truth, are made to recede
from it by time''; the solution he proposes is the use of newspapers:
\begin{quotation}
These daily and periodical publications afford a continual supplement
to the best statistic authors, and form a kind of practical statistics.
\end{quotation}
So Bielfeld's description of the discipline, though quite remote from
the modern meaning, is consistent with the generally accepted
German definitions of the time.
\subsection{The term `statistics'}
That the first use of the word `statistics' in English
appeared in [Bielfeld 1770b], is serendipitous: 
Achenwall's courses could have been imported earlier; 
[B\"usching 1764] could
have been translated instead of [B\"us\-ching 1762]. Other German
works containing the word
were indeed translated between Bielfeld and Sinclair, such as [Zimmermann
1787]. The fact remains that the word had been used in England, more than
twenty years before Sinclair. There is no reason to doubt that
Sinclair discovered it in 1786 ``in the course of a very
extensive tour through the northern parts of Europe'' [Sinclair
1798:xiii]. Yet, [Bielfeld 1770b] had not gone completely
unnoticed. It had been extensively 
reviewed upon appearance in [Smollett 1770:262-270]. More
significantly, the fourth Earl of
Abingdon (1740-1799), in his 1780 \emph{Dedication to the Collective 
Body of the People of England}, after an extensive citation of
[Bielfeld 1770a], added as a footnote:
\begin{quotation}
See Elements of Erudition, vol. 1. p.89, and 103. This Science, in
order to its Attainment, is very judiciously treated of in the third
Volume of these Elements, under the Head of \emph{Statistics}; and
to which the Reader is not only referred in particular, but the
\emph{Elements} themselves in general, as well as \emph{the
Political Institutes} of Baron Bielfeld, are worthy the Perusal of
every Lover of Learning and Science.
[Abingdon 1780:lvii]
\end{quotation}
Even if they did not reach the popularity of the \emph{Institutions
  politiques}, the \emph{Elements of universal erudition} enjoyed a 
considerable success. For instance, it has been shown that they have been
a constant source of inspiration to
Edgar Allan Poe (1809-1849).
\begin{quotation}
The evident popularity of the work in Europe led to a
translation into English by William Hooper, published in London in
1770 in three volumes as Elements of Universal Erudition and pirated
the next year in Dublin. Poe unquestionably used the English
translation for all the items in ``Pinakidia'' except number 154, which
reprints a stanza from a French ``Vaudeville'' 
that Hooper omitted from his version.

Bielfeld's volumes were a major source for many of Poe's learned
allusions, curious bits of information, and even germinal ideas for
tales, and Poe's use of the Empedocles item is important enough, in
Page's treatment, to warrant a close examination of how he drew from
Hooper's translation.  
[Pollin 1980]
\end{quotation}
Writing a vibrant eulogy of his father, Sinclair's son feels quite
uneasy with the question of precedence. He cannot avoid citing Bielfeld.
\begin{quotation}
The vast variety of subjects of which the attention of the statist
should be given, is ably and comprehensively enumerated by Baron
Bielfield in his ``Elements of Universal Erudition''. His work, however,
contains speculations and directions only. He did not attempt to put
his theory in practice by an actual enquiry into the circumstances of
the German empire.
[Sinclair 1837:6]
\end{quotation}
Then he reproduces his father's argument of a semantic shift.
\begin{quotation}
The terms \emph{Statistics}, and \emph{Statistical}, which occurred
continually in this volume, were such novelties in the British
\emph{nomenclature} of economic science, that Sir John thought it
necessary to apologize for their introduction. He explained that he
had derived the term from the German, though he employed it in a sense
somewhat different from its foreign acceptation. In Germany, a
statistical enquiry related to the \emph{political strength} of the
country, or to questions of state policy, whereas he employed the word
to express an enquiry into the state of a country, for the purpose of
ascertaining the amount of \emph{happiness} enjoyed by its
inhabitants, and the means of its future improvement.
[Sinclair 1837:9]
\end{quotation}
Yet, he soon hides away behind the entry ``Statistick'' of Walker's dictionary:
\begin{quotation}
``\emph{Statistick}. This word is not found in any of our Dictionaries. It
seems to have been first used by Sir John Sinclair, in his plan for a
statement of the trade, population, and productions of every parish in
Scotland, with the food, diseases, and longevity of its inhabitants ; a
lan which reflects the greatest credit on the understanding and
benevolence of that gentleman, as it is big with advantages both to
the philosopher and the politician''
[Sinclair 1837:10]
\end{quotation}
He then ackowledges some early remarks about the fact that his
father's use of the word was not essentially different from the German
acception.
\begin{quotation}
German statists, and in particular Professor Schlozer, in his Theorie
der Statistik, insist that a distinction was all along sufficiently
kept in view between politics and statistics, by the statistical
writers of that country, and that Sir John Sinclair's definition was
identical with theirs. I may also here remark, that the Italians may
fairly dispute with their German neighbours the appropriation of this
term, which occurs in some of their writers soon after the revival of
letters.
[Sinclair 1837:10]
\end{quotation}
Let us briefly come back on Sinclair's claim about
changing the meaning of the word,
\begin{quotation}
for by Statistical is meant in Germany, an
    inquiry for the purposes of ascertaining the political strength of
    a country or questions respecting \emph{matters of state}; whereas the
    idea I annex to the term is an inquiry into the state of a
    country, for 
\emph{the purpose of ascertaining the \emph{quantum} of happiness
    enjoyed by its inhabitants, and the means of its future
    improvement}
[Sinclair 1798:xiii]
\end{quotation}
Actually, the `pursuit of happiness' was a recurrent theme in the
philosophy of Enlightment, long before Sinclair. 
As Ian Hacking [1991:194] puts it:
``The fundamental principle of the original moral sciences was
the Benthamite one: the greatest happiness to the greatest
number''. The `happiness of peoples' repeatedly appears 
in the works of Bielfeld, and heads
the first page of the \emph{Institutions politiques}.
\begin{quotation}
Peut-on pr\'etendre que, sans Pr\'eceptes,
les Peuples puissent \^etre constamment heureux dans le cours de
plusieurs Si\'ecles? 
[Bielfeld 1760a:1]
\vskip 1mm
[Can it be claimed that without precepts, peoples can be constantly
happy all along several centuries?]
\end{quotation}
Sinclair's feat lied in the unprecedented scale and accuracy of his
survey; it was rightfully hailed and became an inspiration for the
beginnings of the discipline [Cullen 1975; Hacking 1990]. 
However, his use of the ``new word'' statistic 
is quite far from the colorful ``bare-faced act of robbery''
described by Pearson [1921:2]. Before Pearson, some early historians
of the discipline had given more balanced accounts of the merging of
political arithmetics into statistics. Here is one by 
Walter F. Willcox (1861-1964).
\begin{quotation}
For the historians of statistics are now well agreed that the 
study we cultivate sprang from two main roots, one in the 
German universities, the other in English studies of political 
arithmetic. The former developed, under the name of 
statistics, a descriptive political science almost devoid of figures 
but systematic and suitable for presentation in academic 
lectures or treatises. The latter developed, under the name of 
political arithmetic, a series of fragmentary and disconnected 
studies of available numerical data. Between 1730 and 1830 
the English ideas slowly penetrated Germany, introducing 
numerical data and gaining especially from S\"ussmilch a 
systematic, orderly presentation quite alien to their original 
form. 

During the same period the German name \emph{statistics} spread 
to England and this country. Probably the first writer to 
make it at home in English was Sir John Sinclair whose 
voluminous \emph{Statistical Account of Scotland} exercised a traceable 
influence on both sides of the Atlantic. He wrote in 1798: 
``In the course of a very extensive tour through the northern 
parts of Europe which I happened to take in 1786 I found 
that in Germany they were engaged in a species of political 
inquiry to which they had given the name of Statistics and 
[\ldots] as I thought that a new word might attract more 
public attention I resolved on adopting it and I hope that it 
is now completely naturalized.'' The earliest occurrence of 
\emph{statistics} in English was in 1770 and thus more than fifteen 
years before Sinclair, when Dr. Hooker published a 
translation of Bielfeld's \emph{Elements of Universal Erudition}. One of its 
chapters is entitled Statistics and contains a definition of the 
subject as ``The science that teaches us what is the political 
arrangement of all the modem states of the known world.'' 
With this German name came also some of the German fondness 
for system and for breadth of treatment, and all these 
factors contributed to the establishment of English and 
American statistical societies. The statistics which were thus to 
be studied came much nearer to the German prototype than 
to the English political arithmetic. 
[Willcox 1914]
\end{quotation}
A vivid account of the struggles that accompanied the mutation from
statistik to statistics
is given by John [1883], an article that Pearson had read
(cf. footnote by Pearson's son in [Pearson 1921:3]).
See also [Guy 1865; Hilts 1978; Plackett 1986] 
for etymologies and ancient uses of the words `statist'
and `statism' in English.
 
One question remains to be answered: why did Sinclair's own version
prevail up to this date [Johnson \& Kotz 1993:70-72; Hald 2003:82]?
Sinclair's achievement had an immediate impact and was cited as a
model outside England, in particular in France [Playfair
1802:157]. Beyond the impressive feat, Sinclair's dense network of
correspondents certainly played an important role in the
diffusion. Traces can be found in part VII, volume 1 of Sinclair's 
voluminous correspondence.
For instance, here is a letter of Professor Zimmerman of  Brunswick, 
dated 17th July 1792 (the author of [Zimmerman 1787]). 
\begin{quotation}
I shall not delay a moment to insert an ample extract, in the last
number of my Geographical and Statistical Journal, which I have
published for above two years. These sciences will gain much by your
enterprise; and I feel the greatest anxiety to see a work, of such
extent and utility, brought to a conclusion.
[Sinclair 1831:288]
\end{quotation}
Another letter, from Dr Guthrie, dated at St Petersburg, 
26th september 1792, well describes the early 
reception of the ``statistical
account of Scotland''. 
\begin{quotation}
Your Statistical Work is, in my opinion, the most perfect which has
yet appeared, and will probably serve as a model to other countries,
although few possess the same set of respectable pastors, to collect
materials, living with their flock in habits of friendship and
intimacy, the natural result of the sensible regime of the Scotch
church.
[Sinclair 1831:289]
\end{quotation}
Bielfeld does not seem to have maintained such an important
correspondence, and his relatively short life did not leave him much
time to ensure his own publicity; whereas Sinclair lived long enough to
become the oldest active member (at the age of eighty) of the newly
founded Statistical Society of London in 1834.

\bibliographystyle{decsi}

\end{document}